\def\year{2022}\relax
\documentclass[letterpaper]{article} 
\usepackage{aaai22}  
\usepackage{times}  
\usepackage{helvet}  
\usepackage{courier}  
\usepackage[hyphens]{url}  
\usepackage{graphicx} 
\usepackage{natbib}  
\usepackage{caption} 
\DeclareCaptionStyle{ruled}{labelfont=normalfont,labelsep=colon,strut=off} 
\usepackage{algorithm}
\usepackage{algorithmic}

\usepackage{amsmath}
\usepackage{amsfonts}
\newtheorem{definition}{Definition}

\newtheorem{assumption}{Assumption}

\newtheorem{corollary}{Corollary}

\newtheorem{theorem}{Theorem}

\usepackage{newfloat}
\usepackage{listings}
\floatstyle{ruled}
\newfloat{listing}{tb}{lst}{}
\floatname{listing}{Listing}
\title{Sample Average Approximation for Stochastic Optimization with Dependent Data: Performance Guarantees and Tractability}
\author {
	Yafei Wang,\textsuperscript{\rm 1}
	Bo Pan,\textsuperscript{\rm 1}
	Wei Tu,\textsuperscript{\rm 2}
	Peng Liu,\textsuperscript{\rm 3}
	Bei Jiang,\textsuperscript{\rm 1}
	Chao Gao,\textsuperscript{\rm 4} 
	Wei Lu,\textsuperscript{\rm 5}\\
	Shangling Jui,\textsuperscript{\rm 6}
	Linglong Kong\textsuperscript{\rm 1}\thanks{Corresponding author}
}
\affiliations {
	\textsuperscript{\rm 1}University of Alberta,
	\textsuperscript{\rm 2}Queen's University, 
	\textsuperscript{\rm 3}University of Kent,
	\textsuperscript{\rm 4}Huawei Canada Research Center,\\
		\textsuperscript{\rm 5}Huawei Canada,
	\textsuperscript{\rm 6}Huawei Technologies Ltd.\\
	\{yafei2, pan1, bei1, lkong\}@ualberta.ca, wei.tu@queensu.ca, p.liu@kent.ac.uk,
	chao.gao4@huawei.com,
	robin.luwei@hisilicon.com,
	jui.shangling@huawei.com
}

\usepackage{bibentry}

\begin{document}

\maketitle

\begin{abstract}
Sample average approximation (SAA), a popular method for tractably solving stochastic optimization problems, enjoys strong asymptotic performance guarantees in settings with independent training samples. However, these guarantees are not known to hold generally with dependent samples, such as in online learning with time series data or distributed computing with Markovian training samples. In this paper, we show that SAA remains tractable when the distribution of unknown parameters is only observable through dependent instances and still enjoys asymptotic consistency and finite sample guarantees. Specifically, we provide a rigorous probability error analysis to derive $1 - \beta$ confidence bounds for the out-of-sample performance of SAA estimators and show that these estimators are asymptotically consistent. We then, using monotone operator theory, study the performance of a class of stochastic first-order algorithms trained on a dependent source of data. We show that approximation error for these algorithms is bounded and concentrates around zero, and establish deviation bounds for iterates when the underlying stochastic process is $\phi$-mixing. The algorithms presented can be used to handle numerically inconvenient loss functions such as the sum of a smooth and non-smooth function or of non-smooth functions with constraints. To illustrate the usefulness of our results, we present several stochastic versions of popular algorithms such as stochastic proximal gradient descent (S-PGD), stochastic relaxed Peaceman--Rachford splitting algorithms (S-rPRS), and numerical experiment.
\end{abstract}

\section{Introduction}
Stochastic optimization, a powerful modeling paradigm in optimization under uncertainty, is ubiquitous in statistical machine learning, engineering, and decision-making problems ~\citep{franklin2005elements, heyman2004stochastic, fouskakis2002stochastic}. Specifically, these problems seek to minimize an expected loss taken with respect to the distribution $\mathbb{P}$ of a random parameter $\xi$. However, more often than not, this probability distribution is unknown and can only be observed through a finite number of sample points. We are thus forced to solve a surrogate optimization problem constructed by the observed data: the optimal value of the original problem can be approximated by that of the surrogate problem. The main goal of this paper is to study the properties of sample average approximation (SAA), a powerful approach to stochastic optimization that is considered statistically and computationally ``optimal'' in settings where observations are independent. These claims have the following caveats in non-independent settings.

\begin{itemize}
	\item Most existing statistical guarantees for SAA critically depend on the assumption that training samples are independent and identically distributed (i.i.d.). However, the i.i.d. assumption can be difficult to justify or outright invalid in practice. It is thus important to examine the properties of SAA estimators when samples are known to be correlated.

	\item For an optimization scheme to be useful, it should be solved efficiently. Standard convex optimization techniques, while widely applicable, suffer performance-wise in problems that are complex or highly structured, and trained with non-i.i.d. samples. Consequently,  practically useful methods should offer guarantees that remain valid when the training samples display serial dependencies and be rich enough to handle numerically inconvenient problems.
\end{itemize} 

In this paper, we consider applying SAA to the stochastic optimization problem
\begin{align*} \label{var_obj}\tag{1}
J^{*}&=\underset{x \in \mathbb{X}}{\min}\left\{\mathrm{E}^{\mathbb{P}}[\ell(x;\xi)]=\int_{\Xi} \ell(x;\xi) \mathbb{P}(\mathrm{d} \xi)\right\},
\end{align*}
where $\Xi$ is a sample space with a probability distribution $\mathbb{P}$ and $\mathbb{X}\in \mathbb{R}^d$ is a convex, feasible parameter space. We assume throughout that $\ell(x;\xi)$ is a closed, convex, and proper function and that $\xi\in \Xi$ denotes a sample instance. We let $x^*$ denote the optimal solution to problem \eqref{var_obj}.

In most situations of practical interest, the distribution $\mathbb{P}$ is not known or cannot be efficiently sampled, such as when $\Xi$ is a high-dimensional or combinational sample space~\citep{johansson2007simple,johansson2010randomized}. This restriction removes information essential to solving problem \eqref{var_obj} exactly. We instead consider receiving samples $\{\xi_k\}_{k=1}^K$ from a stochastic process $P=P^k$ indexed by time $k$, where $P$ converges to the stationary distribution $\mathbb{P}$. This is a natural relaxation of the assumption that training samples are i.i.d. following $\mathbb{P}$. As an example, consider $\Xi=\{\xi\in\{0,1\}^d\mid \langle a,\xi\rangle\leq b\}$, where $a\in \mathbb{R}^d, b\in\mathbb{R}$, $\langle a,\xi\rangle = \sum_{i=1}^d a_i\xi_i$, and $\mathbb{P}$ is the uniform distribution over $\Xi$. A straightforward way to obtain a sample from $\mathbb{P}$ is by iterative random sampling from $\{0,1\}^d$ until the constraint on $\Xi$ is satisfied: this approach takes $O(2^d)$ draws to obtain a feasible sample. Alternatively, it is possible to design a Markov chain~\cite{jerrum1996markov} that generates a sample that is $\varepsilon$-close to the distribution $\mathbb{P}$ and only requires $\log(\sqrt{d}/\varepsilon)\exp(O(\sqrt{d}(\log d)^{5/2}))$ draws~\citep{dyer1993mildly}, a greatly reduced sampling cost. Autoregressive processes~\cite{kushner2003stochastic} are yet another example of stochastic data-generating processes but generate a dependent source of data, here the sequential entries of a time series. The assumption of i.i.d. samples is, therefore, unrealistic in many data-generating processes. These two examples highlight the need to consider sampling efficiency and training under inter-sample dependence.

Applications of SAA to stochastic optimization are not new and have been studied extensively in the literature~\citep{kleywegt2002sample,kim2015guide,emelogu2016enhanced,bertsimas2018robust}. As noted above, the idea underlying SAA is simple---to generate solutions to problem \eqref{var_obj}, approximate $\mathbb{P}$ with the discrete empirical distribution $\widehat{\mathbb{P}}_K = \frac{1}{K}\sum_{k=1}^K\delta_{\xi_k}$ corresponding to training samples. SAA improves problem solvability by turning integration over a density function in a summation over discrete points.

Properties of solutions to SAA problems are well understood. In particular, the optimal solution of an SAA problem is known to be strongly consistent and asymptotically normal \cite{kim2015guide}. However, most works study problems in settings where data is i.i.d. 
With the notable exceptions of~\citet{duchi2012ergodic} and \citet{ agarwal2012generalization}, we are not aware of any studies of SAA with non-i.i.d. data, while these work focus more on the convergence of the proposed algorithm and iterate asymptotics rather than the properties of SAA. Results on the asymptotic consistency of SAA under an unknown probability distribution $\mathbb{P}$ and dependent training data have not yet been established and are of particular importance.

Tractability is equally as important as statistical guarantees when establishing the practical utility of an optimization scheme. SAA tractability suffers when the loss function $\ell$ possesses a complex structure, such as statistical machine learning problems that enforce prior knowledge of the form of the solution, such as sparsity, low rank, and smoothness~\citep{franklin2005elements}. As a result, it is critically important to develop algorithms that are both rich enough to capture the complexity of data and scalable enough to process data in a parallelized or fully decentralized fashion. We study the tractability of SAA with non-i.i.d. training samples for a class of stochastic first-order algorithms. We focus primarily on operator-splitting schemes, which are widely used due to their scalability with respect to problem dimensionality. More importantly, operator splitting schemes can be easily parallelized and are usually simple and cheap to implement. 

One general iteration scheme, defined as the Stochastic Krasnosel'ski\u{\i}--Mann (S-KM) iteration, is
\begin{align*}
x^{k} = x^{k-1} + \lambda_{k-1}(T(x^{k-1}) - x^{k-1} + \epsilon_{k-1})
\end{align*}
where $T$ is a nonexpansive operator defined as a mapping:  $T: \mathbb{X} \rightarrow \mathbb{X}$ such that $\lVert T x-T y\rVert \leq \lVert x-y\rVert$ holds for all $x, y \in \mathbb{X}$. The stochastic error $\epsilon_{k-1}$ is caused by uncertainty in random sampling. Succinctly, an operator splitting algorithm converts an optimization problem into a problem of finding a fixed point of a nonexpansive operator and breaks this problem into several relatively simple subproblems. Many commonly used methods, such as stochastic proximal gradient descent (PGD) and  stochastic alternating direction method of multipliers (ADMM), have this iteration step with a specific nonexpansive operator. Classical PGD and ADMM algorithms are special cases of the KM iteration without the stochastic error term.

Our setting is perhaps most similar to that in \citet{duchi2012ergodic}, which also considers receiving data from an ergodic process. However, our work here differs in two fundamental ways. First, \citet{duchi2012ergodic} focuses on algorithmic convergence guarantees with dependent samples but pays little attention to statistical properties of SAA estimators. Second, \citet{duchi2012ergodic} only considers stochastic mirror descent while we consider multiple other algorithms. \citet{sun2018markov} similarly only considers stochastic gradient descent and works under the assumption that training samples are generated by a Markov chain, a special case of an ergodic process. The sampling technique used in \citet{derman2020distributional} is the same as that in \citet{sun2018markov}, except that samples are generated from multiple independent trajectories of serially correlated states. Using multiple replication (MR), a technique that attempts to remove inter-sample dependence via multiple stochastic process trajectories, the authors generate i.i.d. samples and train the model using standard algorithm. This approach 
can help to get rid of dependence among the samples but may be difficult to acquire in practice (we refer to the next section for details). In contrast to this approach, we establish asymptotic consistency and non-asymptotic bounds for SAA problems and study, from a fixed-point iteration perspective, properties of a variety of first-order algorithms based on a \textit{single} trajectory of a stochastic process.

In summary, while many existing works indicate that SAA estimators are asymptotically consistent and tractable, hardly any existing statistical performance guarantees apply when training data fails to be i.i.d. or when training samples are generated from a single trajectory. This paper addresses this gap. Specifically, we study the performance of SAA estimators when training samples are generated by an ergodic stochastic process and, in the same setting, establish properties of the S-KM iteration in solving SAA. The main contributions of this paper are as follows.

\begin{itemize}
	\item \textbf{Asymptotic consistency.} We generalize SAA problems to scenarios where training samples are correlated and prove that the SAA solutions are asymptotically consistent. 
	\item \textbf{Finite sample guarantees.} By introducing a weakened version of $\phi$-mixing, we establish $1 - \beta$ confidence bounds on the out-of-sample performance based on the optimal solution obtained by minimizing an SAA problem.
	\item \textbf{Tractability.} We examine the performance of first-order algorithms for solving SAA problems in an efficient manner via monotone operator theory. We show that the approximation error of the algorithm is bounded and concentrates around zero, and further establish iterate deviation bounds.
\end{itemize}

\section{SAA with Dependent Data}

To motivate the broad applicability of sampling from a stochastic process, in this section, we begin with an example in distributed optimization under a simple peer-to-peer communication scheme~\citep{johansson2007simple}, where the optimization problem evolves according to a finite-state Markov chain. We then move toward a specific sampling technique, called the multiple replication approach, which generates training samples from multiple trajectories to get rid of dependency among the samples. However, its efficiency can not be guaranteed and it wastes a lot of samples. We then propose an iteration procedure that efficiently uses every obtained sample and works when there is only a single trajectory available.

\subsubsection{Peer-to-Peer Optimization}

Suppose that the distribution $\mathbb{P}$ is supported on a set of $n$ points $\{\xi_1,\ldots,\xi_n\}$ and that there are $n$ processors, each with a convex function $\ell(x;\xi_i)$. The objective is to minimize 
$L(x) = \frac{1}{n}\sum_{i=1}^{n}\ell(x;\xi_i)$. To solve this problem, the current set of parameters $x^k\in \mathbb{X}$ is passed to one of the processors and updated in each iteration. More specifically, let the token $i(k)$ indicate the processor $i$ holds $x^k$ at time $k$: at this time, only the data stored in processor $i$ is accessed. Given the current state $i(k)$, the next state $i(k+1)$ is determined randomly via $\mathbb{P}(i(k+1)=j\mid i(k)= i) = P_{ij}$, with $0\leq P_{ij}\leq 1$. The data stored in processor $j$ is then accessed and used to update $x^k$.  Since $\Xi$ is a combinational space in this setting, it is hard to draw samples directly. Moreover, it is unrealistic to assume that samples are i.i.d. However, because the token $i(k)$ can be viewed as evolving according to a Markov chain with a doubly stochastic transition matrix $(P_{ij})_{i,j}$, the data generating process forms a Markov chain.

\subsubsection{Multiple Replication Approach}

As mentioned previously, although we can design a stochastic process to generate training samples, it may be impossible to generate i.i.d. training samples. A natural method called multiple replication approach~\citep{gelman1992inference}, is adopted to obtain a sequence of i.i.d. samples. Specifically, in this approach, after specifying the initial conditions of the stochastic process $P$, a sequence $\xi_1,\ldots,\xi_s$, for some $s$, is generated. And we only keep the last sample $\xi_s$ that follows the marginal distribution $P^s$, which is assumed to be close to $\mathbb{P}$. The same procedure is repeated for $K$ times to simulate $K$ i.i.d. samples, then use standard algorithms as for independent data.  

Unfortunately, this method does not work if there is only one trajectory or expensive to simulate multiple  trajectories. Further, the computation cost could be quite high.  This is because sampling a long trajectory and using only the last sample wastes a large number of samples, especially when $s$ is large. This waste may seem necessary because a small $s$ induces a large bias in $\xi_{s}$: after all, a random trajectory may take a long time to explore the parameter space and will often double back to previously visited states. This further complicates the problem of choosing a $s$ appropriately. A small $s$ will cause large bias in $\xi_{s}$, which slows the convergence of algorithms and reduces its final accuracy---$\{\xi_k\}_{k=1}^K$ are generated from $P^s$ where $P^s$ could be far away from $\mathbb{P}$. A large $s$, on the other hand, is wasteful especially when the iterate $x^k$ is still far from convergence and some bias does not prevent the iteration update to make good progress. Therefore, $s$ should increase adaptively as $k$ increases---this makes the choice of $s$ even more difficult.

\subsubsection{Stochastic Krasnosel'ski\u{\i}--Mann (S-KM) Iteration }

Assume that $\mathbb{P}$ can only be observed through samples $\{\xi_k\}_{k=1}^K$ from an ergodic stochastic process $P$ that converges to $\mathbb{P}$. We address the SAA problem
\begin{align*}\tag{2} \label{obj_emp}
\widehat{J}^{*}_{K} =\underset{x\in \mathbb{X}}{\text{min}}\left\{L(x) = \mathbb{E}^{\widehat{\mathbb{P}}_{K}}[\ell (x; \xi)]=\frac{1}{K} \sum_{k=1}^{K} \ell\left(x; \xi_{k}\right)\right\} 
\end{align*}
and its corresponding optimal solution $\widehat{x}^{*}_K$ with an alternative iteration procedure that uses every sample immediately. Specifically, suppose that the distribution $\mathbb{P}$ is supported on $\{\xi_{k}\}_{k=1}^K$. Problem \eqref{var_obj} can then be approximated by \eqref{obj_emp}. Iteration procedure to solve problem \eqref{obj_emp} is given in Algorithm 1: in the $t$-th iteration, the update
$$x^{k} \leftarrow x^{k-1} + \lambda_{k-1}\big(T(x^{k-1}; \xi_{k}) - x^{k-1}\big)$$
is applied, where $\xi_k$ is sampled from the stochastic process $P$ evaluated at time $k$.The operator $T$ in Algorithm \ref{alg:algorithm} is a nonexpansive operator that depends on the specific method used. We include the sample $\xi_k$ as an argument of $T$ to explicitly indicate that the $k$-th iteration depends only on the most recently drawn sample. 
For more details regarding the forms of $\ell$ and $T$ used in practice, please see the section of Application and Table 1.

\begin{algorithm}[h]
	\caption{Stochastic Krasnosel'ski\u{\i}--Mann (S-KM)}
	\label{alg:algorithm}
	\textbf{Input}: Initial value $x^0$ and given $\delta$-optimality \\
	\textbf{While} \hspace*{-0.1cm}
	$\|T(\bar{x}^{k-1};\xi^{k})-\bar{x}^{k-1}\|^2>\delta$
	\begin{algorithmic}[1] 
		\STATE Sample $\xi^{k} \sim P^k$
		\STATE $x^{k} \leftarrow \bar{x}^{k-1} + \lambda_{k-1}\big(T(\bar{x}^{k-1}; \xi^{k}) - \bar{x}^{k-1} \big)$
		\STATE $\bar{x}^{k} \leftarrow \frac{k - 1}{k} \bar{x}^{k-1} + \frac{1}{k} x^k$ 
	\end{algorithmic}
	\textbf{end while} 
\end{algorithm}

\section{Statistical Guarantees} \label{sg}
\noindent
In this section, we show that the widely used SAA method retains its statistical guarantees if training samples are generated by an ergodic stochastic process that converges to a desired stationary distribution $\mathbb{P}$. The reason for this is that the empirical distribution is a sufficient statistic and satisfies requirements for large number theory.We begin with finite sample performance. If we only have access to $K$ training samples, we can obtain the optimal solution $\widehat{x}_K^{*}$ and the corresponding optimal value $\widehat{J}_K^{*}$ via \eqref{obj_emp}. The quality of $\widehat{x}_K^{*}$ and $\widehat{J}_K^{*}$ can be evaluated through out-of sample performance, defined as $\mathbb{E}^{\mathbb{P}} [\ell (\widehat{x}^{*}_{K};\xi_{\text{test}})]$, where $\xi_{\text{test}}$ is a testing sample assumed to be drawn from $\mathbb{P}$ and is independent of the training samples. 

\subsection{Preliminaries}

We start the section by recalling some definitions that facilitate to present theoretical properties in the next section. The definition of total variation distance is first introduced to measure the convergence of the stochastic process $P$ to the distribution $\mathbb{P}$.

\begin{definition} Let $\mathbb{P}$ and $\mathbb{Q}$ be probability measures defined on a set $\Xi$ with respective densities $p$ and $q$ relative to an underlying measure $\mu$. The total variation distance between $\mathbb{P}$ and $\mathbb{Q}$ is
\begin{equation*}
	d_{\mathrm{TV}}(\mathbb{P}, \mathbb{Q}) = \frac{1}{2}\int_{\Xi}\lvert p(\xi)-q(\xi)\lvert  \mathrm{d} \mu(\xi)= \sup _{A \subset \Xi}\rvert\mathbb{P}(A)-\mathbb{Q}(A)\rvert,
	\end{equation*}
\end{definition}
where the supremum is taken over measurable subsets of $\Xi$.

Using total variation distance, we can define the notion of a mixing stochastic process. Let ${P}^{k}_{[s]}={P}^{k}(\cdot\mid \mathcal{F}_s)$ denote the distribution of $\xi_{k}$ conditional on the $\sigma$-algebra $\mathcal{F}_s$ with $\mathcal{F}_s = \sigma(\xi_1,\ldots,\xi_s)$.

\begin{definition}Define $\mathcal{F}_{0}=\{\emptyset, \Omega\}$ and let $\left(\mathcal{F}_{k}\right)_{k=1}^K$ be an increasing sequence of $\sigma$-algebras such that $\mathcal{F}_{k-1} \subseteq \mathcal{F}_{k}$ for any $k$. The $\phi$-mixing coefficient of the sample distribution $P$ under total variation is
	\begin{equation*}
	\phi(l)=\sup _{k \in\mathbb{N}^+, B \in \mathcal{F}_{k}}\{2 d_{\mathrm{TV}}\big({P}^{k+l}(\cdot\mid B), \mathbb{P}\big)\}.
	\end{equation*}
\end{definition} 
We say that the process is $\phi$-mixing if $\phi(l) \rightarrow 0$ as $l \rightarrow \infty$. Note that if the training samples are i.i.d., then $\phi(1) = 0$. We state the following results in a general form using $\phi$-mixing coefficients.

\subsection{Assumption and Main results}
Before formalizing any statistical properties, we introduce one assumption.
\begin{assumption}
	The $\phi$-mixing coefficients for the sample distribution are summable, i.e.,  $\sum_{k = 1}^{\infty} \phi(k) < \infty$.
\end{assumption}

Assumption 1 is met by some stochastic processes satisfying geometric mixing since $\phi(k)\leq \phi_0\exp(-\phi_1 k^{\alpha})$ would hold for some $\phi_0>0$, $\phi_1>0$, and $\alpha>0$. A large class of stochastic process are geometric mixing: this includes autoregressive models and aperiodic Harris-recurrent Markov processes~\citep{modha1996minimum}.

Let $\lVert\widehat{\mathbb{P}}_K - \mathbb{P}\rVert = \int_0^1\lvert\widehat{\mathbb{P}}_K(t) - \mathbb{P}(t)\rvert\mathrm{d}t$ with $\widehat{\mathbb{P}}_K = \frac{1}{K}\sum_{k=1}^K\delta_{\xi_k}$. If Assumption 1 holds, then Theorem 1 in \citet{dedecker2007empirical} indicates that 
\begin{equation*}
\mathbb{P}\left\{\lVert\widehat{\mathbb{P}}_K - \mathbb{P}\rVert\geq \varepsilon\right\}\leq 2\exp \left(-\frac{K^2\varepsilon^2}{2 C(\sum_{k=1}^{K}\phi(k))}\right) 
\end{equation*}
for all $K\geq 1$ and $\varepsilon>0$, where $C(\sum_{k=1}^{K}\phi(k))$ is a function of $\sum_{k=1}^{K}\phi(k)$ and satisfies $C(\sum_{k=1}^{K}\phi(k)) <\infty$. This concentration inequality provides a prior estimate of the distribution $\mathbb{P}$ that resides outside of the $\varepsilon$-ball $\mathbb{B}_\varepsilon(\widehat{\mathbb{P}}_K) = \{\widetilde{\mathbb{P}}\mid \lVert\widehat{\mathbb{P}}_K - \widetilde{\mathbb{P}}\rVert\leq \varepsilon\}$. Therefore, taking $\varepsilon$ as  
\begin{align}
\varepsilon_K(\beta)= \left(\frac{2C(\sum_{k=1}^{K}\phi(k)) \log(2\beta^{-1})}{K^2}\right)^{\frac{1}{2}}, \label{error} \tag{3}
\end{align}
we get the smallest ball that contains $\mathbb{P}$ with confidence $1-\beta$ for some prescribed $\beta\in(0,1)$.

\begin{theorem}\label{th1} (Out-of-sample guarantees) Let $\widehat{J}^{*}_{K}$ and  $\widehat{x}^{*}_{K}$ be as defined in \eqref{obj_emp}. Suppose that $\ell(x;\xi)$ is bounded by a constant $L$ for $x\in \mathbb{X}$ and $\xi \in \Xi$. Let $\varepsilon^*=L\varepsilon_K(\beta)$. Then, we have that
	\begin{equation*}
	\mathbb{P} \left\{\mathbb{E}^{\mathbb{P}} [\ell (\widehat{x}^{*}_{K};\xi_{\mathrm{test}})] \leq \widehat{J}^{*}_{K} + \varepsilon^*\right\} \geq 1 - \beta.
	\end{equation*}
\end{theorem}

Equation \eqref{error} indicates that $\varepsilon^{*}\rightarrow 0$ as $K\rightarrow\infty$ for any fixed $\beta$. Since the true distribution $\mathbb{P}$ is unknown, the out-of-sample performance of $\widehat{x}^*_K$, defined as $\mathbb{E}^{\mathbb{P}} [\ell (\widehat{x}^{*}_{K}, \xi)]$, cannot be evaluated in practice. It is then more practical to establish bounds on $\mathbb{E}^{\mathbb{P}} [\ell (\widehat{x}^{*}_{K}; \xi)]$. It can be seen directly that $J^*\leq \mathbb{E}^{\mathbb{P}} [\ell (\widehat{x}^{*}_{K};\xi)]$, but this lower bound is still impractical unless $\mathbb{P}$ is known. Our primary concern here is to bound the cost from above. From Theorem \ref{th1}, we can conclude that the out-of-sample performance of $\widehat{x}^*_K$ is bounded by a ball of $\widehat{J}^{*}_{K}$ with radius $\varepsilon^*$ with probability $1-\beta$. \citet{esfahani2018data} establishes a similar results for minimization-maximization problems in the i.i.d. setting. In addition, one can show that if $\beta_K$ converges to zero at a particular rate, then the solution to problem \eqref{obj_emp} converges to the original solution of problem \eqref{var_obj} as $K$ tends to infinity. 

\begin{theorem}(Asymptotic consistency)
	Let $\beta_K \in (0, 1)$ with $\lim_{K \rightarrow \infty} \varepsilon_K(\beta_K) = 0$ and $\sum_{K = 1}^{\infty} \beta_K < \infty$.
	Under the assumptions of Theorem \ref{th1}, we have 
	\begin{equation*}
	\mathbb{P}\left\{\lim _{K \rightarrow \infty} \widehat{J}^{*}_{K}=J^{*}\right\}=1 \quad\text{and}\quad\mathbb{P}\left\{\lim _{K \rightarrow \infty} \widehat{x}^{*}_K=x^{*}\right\}=1.
	\end{equation*}
\end{theorem}

\section{Computational Tractability}
\noindent
Even though SAA offers powerful statistical guarantees, it is practically useless unless the underlying optimization problem can be solved efficiently. In this section, we develop a numerical procedure to solve problem \eqref{obj_emp} when the data comes from an ergodic process $P$ that converges to $\mathbb{P}$. We consider two types of problems related to problem \eqref{obj_emp}: one unconstrained problem, given by
\begin{align}
\underset{x \in \mathbb{X}}{\min}\left\{\mathbb{E}^{\widehat{\mathbb{P}}_{K}}[\ell (x;\xi)]=\frac{1}{K} \sum_{k=1}^{K} f\left(x;\xi_{k}\right) + g\left(x;\xi_{k}\right)\right\}, \label{c1} \tag{4}
\end{align}
and another subject to a linear constraint,
\begin{align} 
&\underset{x \in \mathbb{X}, y \in \mathbb{X}}{\min} \quad  \left\{\mathbb{E}^{\widehat{\mathbb{P}}_{K}}[\ell (x,y;\xi)]=\frac{1}{K} \sum_{k=1}^{K} f\left(x;\xi_{k}\right) + g\left(y;\xi_{k}\right)\right\} \nonumber\\
&\text{subject to  }  A x+B y=b, \label{c2} \tag{5}
\end{align}
where $b \in \mathbb{X}$ and the operators $A$, $B$ are bounded and linear.

Many optimization problems can be cast as one of problem \eqref{c1} or \eqref{c2}~\citep{zhang2004solving,teo2010bundle,davis2019stochastic,yu2019sparse}. Problems of these types arise in diverse applications in image processing, machine learning, and statistics~\citep{boyd2004convex,pietrosanu2020advanced,pietrosanu2021estimation,wang2019wavelet,zhang2021high}. In these fields, the dimensionality of data can be extremely large. Traditional methods may thus fail to efficiently (in terms of time) generate solutions. Regularizers, for example, enforce prior knowledge of the form of the solution, such as sparsity, low rank, or smoothness. In regularization schemes, $f$ and $g$ can be data fitting and penalty terms, respectively.  Typically, penalty terms make problems \eqref{c1} and \eqref{c2} difficult to optimize jointly. Even if the both terms can be handled jointly, modern data is often high-dimensional and consists of millions or billions of training examples: running even a single iteration using classical algorithms is often infeasible. Moreover, in most statistical learning problems, we are more concerned with target parameter estimates rather than the objective function value. We present a stochastic KM algorithm that can handle a amount of non-i.i.d. data in a fast, parallelized, and efficient manner.

The methodology we present is different from those used in classical convex analysis~\citep{boyd2004convex}, mainly because operator splitting algorithms are driven by fixed-point iteration rather than by the goal to minimize a loss function: convergence is due to the contraction property of a given fixed-point operator instead of ``descent'' on the loss. Most fixed-point iteration schemes do not decrease the objective function monotonically. Therefore, convergence of the objective function is a consequence of fixed-point convergence but not the cause of it.

For a nonexpansive operator $T: \mathbb{X} \rightarrow \mathbb{X}$, define $\operatorname{Fix}(T)=\{x\in \mathbb{X}: x = T(x)\}$. We assume that $\operatorname{Fix}(T) \neq \emptyset$. Let $\lambda_{k} \in(0,1)$ and choose $x^{0}$ arbitrarily from $\mathbb{X}$. Then the S-KM iteration of $T$ with data generated from $P$ at time $k$ is
\begin{align*} 
x^{k}& = x^{k-1}+\lambda_{k-1}\left(T\left(x^{k-1}\right)-x^{k-1} + \epsilon_{k-1}\right) \\
&= T_{\lambda_{k-1}} (x^{k-1}) + \lambda_{k-1}\epsilon_{k-1},
\end{align*}
where $\epsilon_{k-1}$ is caused by the randomness of samples.

The convergence of the fixed point iteration with a nonexpansive operator $T$ fails in general. The S-KM algorithm thus replaces $T$ with an averaged version $T_{\lambda}$ to ensure convergence, as an averaged nonexpansive operator has the contraction property~\citep{davis2016convergence}. It can be shown that the fixed point of a nonexpansive operator is also the fixed point of its corresponded averaged nonexpansive operator. In practice, the operator $T$ in Algorithm 1 depends on the stochastic splitting method used. For example, when using stochastic proximal gradient algorithm to solve problem \eqref{c1}, $T = \mathcal{J}_{\gamma \partial g}\circ (I - \gamma\partial f)$, where $\mathcal{J}_{\gamma \partial g} = (I + \gamma\partial g)^{-1}$, $I$ is an identity operator, $\gamma$ is a step size. When $g$ is a closed, convex, and proper function, $\mathcal{J}_{\gamma \partial g} $ is equivalent to the well-known proximal operator
$$\operatorname{prox}_{\gamma g}(x)=\underset{y \in \mathbb{X}}{\arg \min } (g(y)+\frac{1}{2 \gamma}\lVert y-x\rVert^{2}_2),$$ with $\lVert\cdot\rVert_2$ denoting the $L_2$ norm on $\mathbb{X}$, and, the corresponding algorithm is known as proximal point approach (PPA). More examples are deferred to Table 1.

\begin{table*}[h] 
	\centering \scriptsize
	\caption{Overview of several first-order algorithms}
	\begin{tabular}{ccc}
		\hline
		Algorithm & Operator identity &  Subgradient identity  \\ 	\hline
		SGD ($g=0$) & $I - \gamma \nabla f$ & $x^{k+1} = x^{k} - \gamma_k \nabla f(x^{k})$  \\
		PPA ($g=0$)  & $(I + \gamma \partial f)^{-1}$ & $x^{k+1} = \operatorname{prox}_{\gamma_{k} f} (x^{k})$\\
		PGD & $(I + \gamma \partial g)^{-1}(I - \gamma \nabla f)$ & $x^{k+1} = \operatorname{prox}_{\gamma_{k} g} (I - \gamma_{k} \nabla f (x^{k}))$ \\
		DRS & $\left(I+\gamma  \partial f\right)^{-1}\left[\left(I+\gamma \partial g\right)^{-1}\left(I-\gamma \partial f\right)+\gamma \partial f\right]$ & $x^{k+1}=\frac{1}{2} x^{k}+\frac{1}{2} \operatorname{refl}_{\gamma_k \partial f} \circ \operatorname{refl}_{\gamma_k \partial g}\left(x^{k}\right)$  \\
		Relaxed PRS & $\left(I+\gamma \partial f\right)^{-1}\left(I-\lambda \partial g\right)\left(I+\lambda \partial g\right)^{-1}\left(I-\lambda \partial f\right)$ & $x^{k+1}=\left(1-\lambda_{k}\right) x^{k}+\lambda_{k} \operatorname{refl}_{\gamma_k \partial f} \circ \operatorname{refl}_{\gamma_k \partial g}\left(x^{k}\right)$ \\ 	\hline
	\end{tabular}
\end{table*}

\section{S-KM Performance with Non-i.i.d. Data}

We next study the properties of S-KM iteration when using non-i.i.d. training samples. We show that, under some mild conditions, S-KM iterates concentrate around the true value. These general results are fundamental and cover many splitting algorithms as special cases. Because splitting algorithms are driven by fixed-point operators, it becomes natural to perform the analysis in terms of Fixed Point Residual (FPR)~\citep{davis2016convergence}, defined as 
$$
e_{k}^2=\lVert Tx^{k}-x^{k}\rVert^2,
$$
which are related to differences between successive KM iterates through $x^{k+1} - x^{k} = \lambda_k(T(x^k) - x^k)$. In first-order algorithms, with the assumption that, $\epsilon_k = 0$, $\forall k$, FPR typically relates to the  gradient of the objective. For example, in the unite-step gradient descent algorithms $x^{k}=x^{k-1}- \nabla f\left(x^{k-1}\right)$, and so the FPR is given by $\lVert\nabla f\left(x^{k-1}\right)\rVert^{2}$. Thus, FPR convergence naturally implies the convergence of $\lVert x^{k+1} - x^{*}\rVert^2$. 

We proceed by establishing the properties of the ergodic FPR. We first, through Theorem \ref{th3}, provide the following boundedness on the expectation of the norm of the approximation error due to randomness of sampling from $P$ rather than $\mathbb{P}$ .

\begin{assumption}
	$\mathbb{X}$ is compact and has finite radius $r$: specifically, for any $x, x^*\in \mathbb{X}$, $\lVert x-x^*\rVert\leq r<\infty$. 
\end{assumption} 

Assumption 2 is the same as the one for online algorithms with correlated data~\citep{agarwal2012generalization,sun2018markov} and common in the online learning, optimization literature. 

\begin{theorem} (Boundedness on approximation error)\label{th3}
	Under Assumption 2, the norm of the difference between the true function $T(x^k;\xi) - x^k$ with $\xi$ drawing from $\mathbb{P}$ and its approximation $T(x^k;\xi_{k+1}) - x^k$ with $\xi_{k+1}$ drawing from $\widehat{\mathbb{P}}_K$ is uniformly bounded in expectation. Specifically, $\mathbb{E}\lVert \epsilon_{k}\rVert \leq \Delta$, where
	\begin{align*}
	\Delta = \left(\frac{8r^2C(\sum_{k=1}^{\infty}\phi(k))}{K^2}\right)^{\frac{1}{2}} \Gamma\left(\frac{1}{2}\right)
	\end{align*}
	and $\Gamma(z) = \int_{0}^\infty x^{z-1}\exp(-x)\mathrm{d}x$ is the gamma function.
\end{theorem}

Theorem \ref{th3} suggests that the approximation becomes close to the true one when $K$ increases. However, the noisy can be large when $K$ is small. This is because we are considering drifting distributions and, in the worse case, $P^k$ can be quit far away from $\mathbb{P}$. Therefore, both Theorem \ref{th3} and MR approach indicate that underestimating the mixing time can potentially backfires.

\begin{theorem} (Bound of fixed point residual) \label{th4}
	Let $\bar{e}_{K}={\Lambda_{K}}^{-1}\sum_{k=1}^{K} \lambda_{k} e_{k}$, with $e_{k}=T x^{k}-x^{k}$, $\Lambda_{K}=\sum_{k=1}^{K} \lambda_{k}$, and $\lambda_{k} \in(0,1)$. Under Assumption 2,
	\begin{equation*}
	\mathbb{E}\lVert\bar{e}_K\rVert\leq \frac{2r(1+\phi(1))+ 2\sum_{k=1}^{K}\mathbb{E}\lVert\lambda_k\epsilon_k\rVert}{\Lambda_K}.
	\end{equation*}
\end{theorem}

When the data is i.i.d. following the distribution $\mathbb{P}$, then $\phi(1)=0$. The result in Theorem \ref{th4} consequently reduces to that for S-KM iterations with i.i.d. samples. To establish an upper bound for S-KM iterates around the true value, we introduce the following two assumptions.

\begin{assumption}
	There is a non-increasing sequence $\kappa(k)$ such that, if $x^k$ and $x^{k+1}$ are successive S-KM iterates, then $\mathbb{E}\{\lVert x^{k+1}-x^k\rVert\mid\mathcal{F}_t\}\leq \kappa(k)$.
\end{assumption}

\begin{assumption}
	For a sequence of samples $\xi_1,\ldots,\xi_K$, the S-KM iteration produces a sequence of iterates $x^1,\ldots,x^{K-1}$ such that $\sum_{k=0}^{K-1}\lVert T_{\lambda_{k}}(x^k;\xi_{k+1}) - T_{\lambda_{k}}(x^*;\xi_{k+1})\rVert\leq R_{K-1}$.
\end{assumption}

Assumption 3 ensures that S-KM iterates are approximately stable. A similar condition for the non-i.i.d. setting is also given in~\citet{agarwal2012generalization}. Assumption 4 allows us to quantify the impact of our assumptions on the performance of some specific instances of $T$.

\begin{theorem} \label{th5} (Deviation bound for iterates)
	Under Assumptions 2, 3, and 4, for any $\tau>0$,
	\begin{align*}
	\mathbb{E}&[\lVert\sum_{k=1}^{K}(x^k - x^*)\rVert] \leq (1+\phi(1))\mathbb{E}[R_{K-1}] \\
	&+ 2(K-\tau) r \sqrt{\phi(\tau+1)} + \tau (\sum_{k=1}^{K-\tau}\mathbb{E}[\kappa(k-1)] + r).
	\end{align*}
\end{theorem}
In the case where $\tau = 0$, $\phi(1) = 0$ gives a bound for the i.i.d. setting.

\section{Application} \label{subsec:appl}

There are many works that focus on using stochastic operator splitting algorithms to solve structured optimization problems~\citep{xu2020primal, rosasco2019convergence,yun2020general, ouyang2013stochastic}. We next present several examples of nonexpansive operators that cover widely used algorithms based on the computation of proximal and gradient operators. For simplicity, we assume that step size $\gamma = \gamma_k$ for all $k$.

\subsubsection{Stochastic PGD} 

Suppose that the function $f$ in problem \eqref{c1} is convex and differentiable with a $(1/\beta)$-Lipschitz continuous gradient for some $\beta > 0$ and that $g:\mathbb{X} \rightarrow \mathbb{R} \cup \{\infty\}$ is a  proper, closed, lower semi-continuous convex function. Solving problem \eqref{c1} is equivalent to finding $x \in \mathbb{X}$ such that $0 \in \partial f (x) + \partial g (x)$. Stochastic PGD (S-PGD), due to its simplicity, efficiency, and empirical performance, is commonly used to solve this problem. For all $k \geq 0$ and $\gamma \in (0, 2\beta)$, the iteration step at time $k$ can be written as
\begin{equation*}
x^{k+1} = \mathrm{prox}_{\gamma \partial g} (x^{k} - \gamma \nabla f(x^{k};\xi_{k+1}) + \epsilon_{f, k})+ \epsilon_{g, k}.
\end{equation*}

We will show that the S-PGD algorithm is a special case of the S-KM iteration. Let $T_{1}=\mathcal{J}_{\gamma \partial g}$ and $T_{2}= (I-\gamma \nabla f)$. Then $T_{\text{PGD}} = T_1\circ T_2$.
Since $T_1$ is (1/2)-averaged and $T_2$ is ($\gamma/(2\beta)$)-averaged, it follows that $T_{\text{PGD}}$ is $2\beta/(4\beta-\gamma))$-averaged~\citep{bauschke2011convex}. Since $\nabla f$ is single-valued, we have that, for all $k\in \mathbb{N}$ and $x\in \mathbb{X}$,
\begin{align*}
x \in(\nabla f+\partial g)^{-1}(0) &\Leftrightarrow x-\gamma \nabla f (x) 
\in x+\gamma \partial g(x)\\
&\Leftrightarrow x \in \operatorname{Fix}(T_{1}\circ T_{2}).
\end{align*}
The results of Theorem \ref{th3}-\ref{th5} then hold. The following corollary further establishes a generalized error bound for regret for S-PGD.

\begin{corollary} (Generalized error bound for regret: S-PGD) Under Assumption 2 and let $\bar{x} = {K}^{-1}\sum_{k=1}^{K} x^k$ and $\underline{\tau}=\inf_{k}\tau_k$, with $\tau_k=\lambda_k(1-\lambda_k)$, we have
	\begin{align*}
	\mathbb{E}&\left\{ f\left(\bar{x}\right)+g\left(\bar{x}\right)-  \left[f(x^{*})+g(x^{*})\right] \right\}
	\\&\leq \frac{r^2}{2K\gamma}+\left(\frac{1}{\beta} - \frac{1}{\gamma}\right)\left(4r^2 + \frac{8r^2\pi C(\sum_{k=1}^{\infty}\phi(k))}{K\underline{\tau}}\right).
   \end{align*} 	
\end{corollary}

\subsubsection{Stochastic generalized DRS}

A line search can be used to guarantee the convergence of the S-PGD algorithm if the Lipschitz constant of $\nabla f$ is not known. Finding an appropriate step size, however, presents another expensive practical challenge. We introduce the Stochastic generalized Douglas–Rachford Splitting (S-gDRS) algorithm to avoid choosing the step size altogether. The results given in Theorem \ref{th3}--\ref{th5} hold by specifying $T$ in the S-KM iteration as
$$
T_{\text{DRS}} = \frac{1}{2}(\mathrm{refl}_{\gamma \partial f}\circ \mathrm{refl}_{\gamma \partial g}+I),
$$
where $\operatorname{refl}_{\gamma\partial f}$ and $\operatorname{refl}_{\gamma \partial g}$ are reflection operators.

\begin{definition} \label{ref}
	Given any operator $T: \mathbb{X}\rightarrow \mathbb{X}$, let $\mathcal{J}_{\gamma T}$ denote the operator $(I+\gamma T)^{-1}$. The operator $\operatorname{refl}_{\gamma T} = 2\mathcal{J}_{\gamma T}-I$ is called the reflection operator of $T$.
\end{definition}

\begin{figure*}[t]
	\centering
	\includegraphics[width=0.9\textwidth]{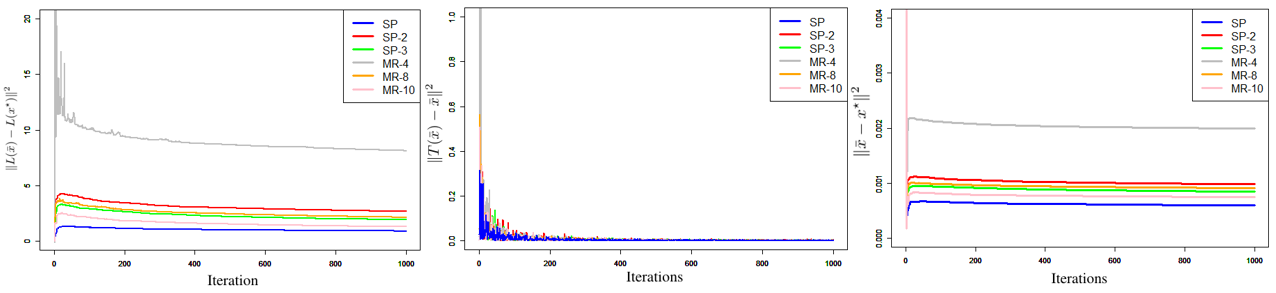} 
	\caption{Performance of Algorithm 1 and the other two methods: SP-m and MR-s.} 
	\label{f1}
\end{figure*}

Because reflection operators are nonexpansive and the composition of nonexpansive operators is nonexpansive~\citep{bauschke2011convex},
we have that $\operatorname{refl}_{\gamma \partial f}\circ \operatorname{refl}_{\gamma \partial g}$ is a nonexpansive operator. Therefore, $T_{\textrm{DRS}}$ is ($1/2$)-averaged, indicating that the S-gDRS algorithm is a special case of the S-KM algorithm. In addition, the definition of the reflection operator indicates that the step size $\gamma>0$ will not affect the convergence of $T_{\text{DRS}}$. This is an advantage over $T_{\text{PGD}}$, which instead requires that $\gamma \in (0,2\beta)$. The S-gDRS algorithm (with $\lambda_k = 1, \forall k$) is a special case of the relaxed Peaceman--Rachford Splitting (PRS) algorithm with the iteration,
$$
x^{k+1}=\left(1-\lambda_{k}\right) x^{k}+\lambda_{k} \cdot \operatorname{refl}_{\gamma \partial f} \circ \operatorname{refl}_{\gamma \partial g}\left(x^{k}\right)
$$
taking $\lambda_k = 1 / 2, \forall k$. A similar result holds for stochastic relaxed PRS (S-rPRS). We also give the error bound for regret for S-rPRS.

\begin{corollary} (Generalized error bound for regret: S-rPRS) Under Assumption 2, Let $\bar{x}^f = K^{-1}\sum_{k = 1}^{K}x_k^f$ and $\bar{x}^g = K^{-1}\sum_{k = 1}^{K}x_k^g$ with the auxiliary points $x_k^{g}=\mathrm{prox}_{\gamma g}(x^k)+\epsilon_g$ and $x_k^{f}=\mathrm{prox}_{\gamma f}\left(\mathrm{refl}_{\gamma \partial g}(x^k)\right)+\epsilon_f$. Denote $\underline{\tau}=\inf_{k}\tau_k$, with $\tau_k=\lambda_k(1-\lambda_k)$.  We have that, 
	\begin{align*}
	\mathbb{E}& [f(\bar{x}^f;\xi) + g(\bar{x}^g;\xi)- (f(x^{*};\xi)+g(x^{*};\xi))]\\
	\leq&\frac{r^2}{4\gamma\lambda K}+ \frac{2(\lambda-1)r^2}{\gamma\lambda^2} 
	+ \frac{4r^2\pi}{\gamma\lambda\underline{\tau}}\left(1-\frac{1}{\lambda}\right)\frac{C(\sum_{k=1}^{\infty}\phi(k))}{K}.\\
	\end{align*}
\end{corollary}

\subsubsection{Numerical Experiments}

We present some numerical results in our setting where samples are generated from an autoregressive process. We show that the iteration procedure in Algorithm 1 uses fewer samples and yields better performance than does the multiple replication approach. We also demonstrate the advantage of using each sample of one trajectory in each iteration rather than at regular intervals: although the dependency of data is weakened by using samples at regular intervals, the performance of iteration has not been improved. Our data-generating mechanism resembles that in \citet{duchi2012ergodic}. Let $A$ be a subdiagonal matrix with entries $A_{i, i-1}\overset{\text{i.i.d.}}{\sim}\mathrm{U}[0.8,0.99]$. We uniformly draw a sparse vector $x \in \mathbb{R}^{1000}$, specifically, with the first non-zero 50 elements of $x$. The data $\{(\xi_{k}^{1}, \xi_{k}^{2})\}_{k\in\mathbb{N}}$ is generated according to the autoregressive process
\begin{align*}
	\xi_{k}^{1}=A \xi_{k-1}^{1}+e_{1} W_{k}, \quad
	\xi_{k}^{2}=\left\langle x, \xi_{k}^{1}\right\rangle+E_{k},\label{ar}
\end{align*}
where $e_{1}$ is the first standard basis vector, the $W_{k}$s are i.i.d. $\mathrm{N}(0,1)$ random variables, and the $E_{k}$s are i.i.d. biexponential random variables with variances of one. We aim to solve the lasso-type problem
\begin{align*} 
	\widehat{x}^{*}_K = \underset{x\in\mathbb{R}^d}{\textrm{argmin}}\left\{\frac{1}{K} \sum_{k=1}^{K}\lVert\left\langle x, \xi^{1}_k\right\rangle-\xi^{2}_k\rVert^2 + \lambda \lVert x\rVert_1\right\},
\end{align*}
where $\lambda$ is a pre-set tuning parameter, using PGD. We generate samples from the above autoregressive model in three ways: (SP) samples are the elements of a single trajectory and are used immediately; (SP-$m$) samples are generated at every $m$-th element of the same trajectory up to get $K$ samples (e.g., with $K = 4$ and $m = 3$, we keep samples at $k=1,4,7,10$); and, (MR-$s$) samples are generated via multiple replication method as the $s$-th elements of independent trajectories starting from the same state. And, we need simulate $K$ trajectories in total. By not using every element, SP-$m$ weakens dependencies between generated samples. SP and MR  are closely resemble the sampling technique in \citet{duchi2012ergodic} and \citet{sun2018markov}. Given sample size, $K = 1000$, we consider $m = 2, 3$ for SP-$m$ and $s = 4, 6, 8, 10$ for MR-$s$.

Figure 1 illustrates our numerical results and the convergence behavior of the three methods that is evaluated by three criteria: regret function $L(\bar{x})-L(x^*)$, FPR and the difference between iterate and true value. As expected, the multiple replication approach shows poor performance for small $s$ as the true mixing time is underestimated: MR-4 has the worst performance under the criteria. Moreover, it becomes clear that using each sample sequentially (SP) rather than attempting to draw weak dependent samples at each iteration from the same trajectory (SP-m) is a more computationally efficient approach.

\section{Conclusion}

In this paper, We show that SAA retains its asymptotic consistency and out-of-sample performance when data is not independent and can be solved efficiently in practice, we also evaluate the performance of a class of first-order algorithms and give several examples illustrating the usefulness of our analyses. We provide generalized error bounds for iterates around the true value that show the impact of dependence of the training sample on the convergence result. It may be possible to sharpen these results using monotone operator properties, such as the contraction property of firmly-nonexpansive operators. We leave these investigations to future work.

\section{Acknowledgements}
	We would like to thank the anonymous reviewers for great feedback on the paper.  Dr. Jiang and and Dr. Kong were supported by the Natural Sciences and Engineering Research Council of Canada (NSERC). Dr. Kong was also supported by the University of Alberta/Huawei Joint Innovation Collaboration, Huawei Technologies Canada Co., Ltd., and Canada Research Chair in Statistical Learning. 

\bibliography{aaai2021}

\end{document}